\documentclass[12pt,dvipdfm]{article}
\usepackage[notcite, notref, final]{showkeys}
\usepackage{amssymb,latexsym}
\usepackage[thmmarks]{ntheorem}
\newtheorem{theorem}{Theorem}[section]

\newtheorem{proposition}[theorem]{Proposition}
\newtheorem{corollary}[theorem]{Corollary}

\newtheorem{conjecture}[theorem]{Conjecture}
\theoremstyle{nonumberplain}
\theoremheaderfont{\it}
\theorembodyfont{\normalfont}
\theoremsymbol{\ensuremath{\Box}}
\newtheorem{proof}{Proof.}
\usepackage{amsmath,amsfonts,latexsym,epsfig}
\usepackage{graphicx, color}
\usepackage[dvips,bookmarks,colorlinks]{hyperref}
\newcommand{\al}{\alpha}
\newcommand{\be}{\beta}
\newcommand{\ga}{\gamma}
\newcommand{\lam}{\lambda}

\newcommand{\om}{\omega}
\newcommand{\vt}{\vartheta}

\newcommand{\rar}{\rightarrow}
\newcommand{\mpt}{\mapsto}
\newcommand{\f}{\frac}
\newcommand{\tf}{\tfrac}
\newcommand{\ld}{\ldots}
\newcommand{\mh}{{\hat{m}}}
\newcommand{\mcC}{\mathcal{C}}

\newcommand{\mcP}{\mathcal{P}}
\newcommand{\mfS}{\mathfrak{S}}
\newcommand{\chihat}{\widehat{\chi}}
\newcommand{\chitilde}{\widetilde{\chi}}
\newcommand{\Aut}{\mathrm{Aut}}
\title{An explicit form for Kerov's character polynomials}
\author{I.P. Goulden\footnote{Department of Combinatorics and
Optimization, University of Waterloo, email: {\small\texttt
ipgoulden@math.uwaterloo.ca}}~
and A. Rattan\footnote{Department of Combinatorics and
Optimization, University of Waterloo, email: {\small\texttt
arattan@math.uwaterloo.ca}}}
\date{April 20, 2005}
\begin{document}
\maketitle

\begin{abstract}
\noindent
Kerov considered the normalized characters of irreducible
representations of the symmetric group, evaluated on a cycle,
as a polynomial in free cumulants. Biane has proved that this
polynomial has integer coefficients, and made various
conjectures.  Recently, \'Sniady has proved Biane's conjectured
explicit form for the
first family of nontrivial terms in this polynomial.
In this paper, we give an explicit expression for all terms in
Kerov's character polynomials. Our method is through
Lagrange inversion.
\end{abstract}

\section{Introduction}
\subsection{Background and notation}
A {\em partition} is a weakly ordered list of positive
integers $\lam =\lam_1\lam_2\ld\lam_k$, where $\lam_1\geq\lam_2\geq\ld\geq\lam_k$.
The integers $\lam_1,\ld ,\lam_k$ are called the {\em parts} of
the partition $\lam$, and we denote the number of parts by $l(\lam)=k$.
If $\lam_1+\ld +\lam_k=d$, then $\lam$ is a partition of $d$, and
we write $\lam\vdash d$. We denote by $\mcP$ the set of all
partitions, including the single partition of $0$ (which
has no  parts).  For partitions $\om, \lam \vdash n$ let $\chi_\om(\lam)$ be
the character of the irreducible representation of the symmetric
group $\mfS_n$ indexed by $\om$, and evaluated on the conjugacy class
 $\mcC_{\lam}$ of $\mfS_n$, which consists of all permutations whose disjoint cycle lengths
are specified by the parts of $\lam$. 

Various scalings of irreducible symmetric group characters have been considered in
the recent literature.  The {\em central character} is given by
\begin{equation*}
\chitilde_{\om}(\lam)=\vert\mcC_{\lam}\vert\f{\chi_\om(\lam)}{{\chi_\om(1^n)}},
\end{equation*}
where $\chi_\om(1^n)$ is
the {\em degree} of the irreducible representation
indexed by $\om$.
For results about the central character, see, for example,~\cite{cgs,fjr,ka}.
Related to this scaling, for the conjugacy class $\mcC_{k1^{n-k}}$ only,
is the {\em normalized character}, given by
\begin{equation*}
  \label{eq:chihatdef}
  \chihat_\om(k1^{n-k}) = n(n-1) \cdots (n-k+1)
\frac{ \chi_\om(k1^{n-k})}{\chi_\om(1^n)}=k\chitilde_{\om}(k1^{n-k}).
\end{equation*}

The subject of this paper is a particular polynomial expression for the
normalized character. The statement of this expression requires some notation involving the
partition $\om$ of $n$. We adapt the following description from
Biane~\cite{bi1,bi2}: consider the Young diagram of $\om$, in the French
convention (see \cite[footnote page 2]{ma}),
and translate it, if necessary, so that the
bottom left of the diagram is placed at the origin of an $(x,y)$ plane.
Finally, 
rotate the diagram counter-clockwise by $45^\circ$.  Note that $\om$ is uniquely
\begin{figure}[ht]\label{fig:secondpart}
 \centering
  \includegraphics[width=12.0cm]{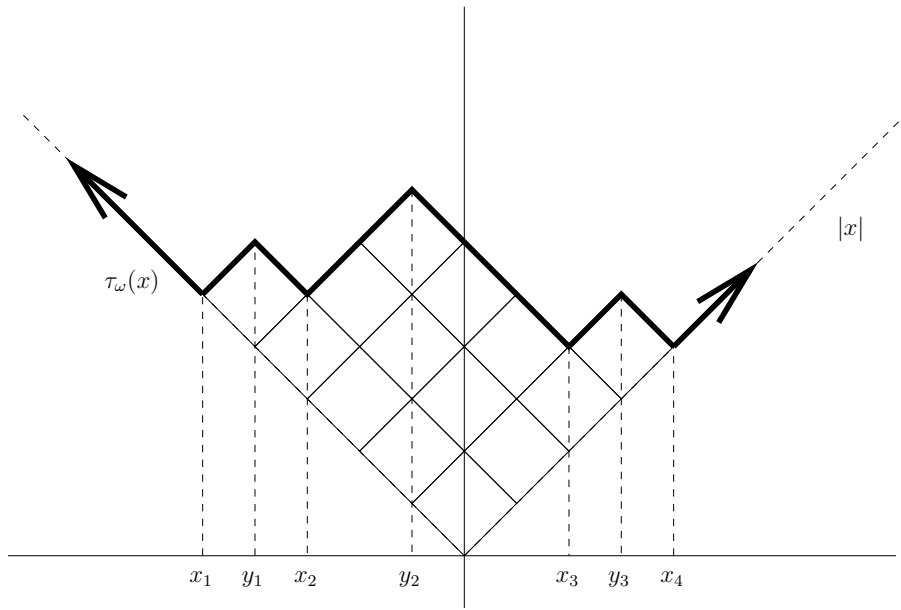}
  \caption{The partition (4\,3\,3\,3\,1) of 14, drawn in the French convention,
  and rotated by $45^\circ$.}
\end{figure}
determined by the curve $\tau_\om(x)$ (see
Figure \ref{fig:secondpart}).  The value of
$\tau_\om(x)$ is equal to $|x|$ for large negative or positive values of $x$ and it
is clear that $\tau_\om^\prime (x) = \pm 1$, where differentiable. The
points $x_i$ and $y_i$ are the
$x$-coordinates of the local minima and maxima, respectively, of the
curve $\tau_\om(x)$.  We suitably scale the size of the boxes in our Young
diagram so that the points $x_i$ and $y_i$ are integers.
Setting $\sigma_\om(x) = (\tau_\om(x) - |x|)/2$, consider the function
\begin{equation}
  \label{eq:gsublam}
  H_\om(z) =\f{1}{z} \exp \int_\mathbb{R} \frac{1}{x-z}
  \; \sigma_\om^\prime(x) \; dx .
\end{equation}
Carrying out the above integration one obtains
\begin{equation*}
H_\om(z) = \frac{\prod_{i=1}^{m-1} (z - y_i)}{\prod_{i=1}^m (z - x_i)} ,
\end{equation*}
where $m$ is the number of nonempty rows in the Young diagram of $\om$.
Now let $R_i(\om)$, $i\geq 1$ be defined by 
\begin{equation*}
z^{-1}+\sum_{i\geq 1}R_i(\om)z^{i-1} = H_{\om}^{\langle -1 \rangle}(z),
\end{equation*}
where $\langle -1\rangle$ denotes compositional inverse.
he $R_i(\om)$'s are known as {\em free cumulants}
in free probability theory.  In this context, they appear in
the asymptotic evaluation of characters.  Specifically, if $\sigma_n
\in \mfS_n$, $n \geq 1$, is a sequence of permutations (subject to a
certain ``balanced'' restriction on the associated Young diagram) with
$k_i$ cycles of length $i$ for $i \geq 2$ and $r = \sum_i
ik_i$, then we have
\begin{equation*}
  \lim_{n \rar \infty} \frac{\chi_\om(\sigma_n)}{\chi_\om( 1^n)} = \prod_{i \geq 2}
  R_{i+1}^{k_i}(\om)n^{-r} + O(n^{-\frac{r+1}{2}}) .
\end{equation*}
For more information about the asymptotics of characters of the
symmetric group (and free cumulants) see, for example, \cite{bi1,io,ke}.

\subsection{Kerov's character polynomials}
The particular polynomials that are the subject of this paper
involve the $R_i(\om)$'s. They first appeared in Biane~\cite{bi2},
where the  following result is stated (as Theorem 5.1).
\begin{theorem}\label{kebi}
For $k\geq 1$, there exist universal polynomials $\Sigma_k$, with
  integer coefficients, such that
  \begin{equation}
    \label{eq:kerovpolys}
    \chihat_\omega(k1^{n-k}) = \Sigma_k(R_2(\omega), R_3(\omega), \ldots,
    R_{k+1}(\om)),
  \end{equation}
for all $\om\vdash n$ with $n\geq k$.
\end{theorem}
Biane attributes Theorem~\ref{kebi} to Kerov, who described this result
in a talk at an IHP conference in 2000, but a proof first appears
in a later paper of Biane~\cite{bi}. The polynomials $\Sigma_k$
are known as {\em Kerov's character polynomials}.
They are referred to as ``universal polynomials'' in Theorem~\ref{kebi}
to emphasize that they are independent of $\om$ and $n$, subject
only to $n\geq k$. Thus we write them with $R_i(\om)$ replaced by
an indeterminate $R_i$, $i\geq 2$.
In indeterminates $R_i$,
the first six of Kerov's character polynomials, as listed in~\cite{bi2}, are
given below:
\begin{align*}
  \Sigma_1 & = R_2\\
  \Sigma_2 & = R_3\\
  \Sigma_3 & = R_4 + R_2\\
  \Sigma_4 & = R_5 + 5R_3\\
  \Sigma_5 & = R_6 + 15R_4 + 5R_2^2 + 8R_2\\
  \Sigma_6 & = R_7 + 35R_5 + 35R_3 R_2 + 84R_3
\end{align*}
Note that all coefficients appearing in this list are positive.
It is conjectured that this holds in general: that for any $k\geq 1$, all nonzero
coefficients in $\Sigma_k$  are positive.

In Biane~\cite{bi2}, this conjecture, which we shall refer to as
the {\em R-positivity conjecture}, is
attributed to Kerov. It has been verified
for all $k$ up to $15$ by Biane~\cite{bi}, who
computed $\Sigma_k$ for $k\leq 15$,
using an implicit formula for $\Sigma_k$ (Theorem 5.1) that he credits
to Okounkov (private communication). Biane further
comments that ``It seems plausible that S. Kerov was aware
of this (see especially the account of Kerov's central limit
theorem in~\cite{io}).'' The following result gives an
adaptation of Biane's formula that appears in Stanley~\cite{st2}.

\begin{theorem}\label{bist}
Let $R(x) = 1 + \sum_{i \geq 2} R_i x^i$ and
\begin{equation}\label{FGR}
F(x)=\f{x}{R(x)},\;\;\;\;\;\;\;\;\;\;
 G(x)=\f{1}{F^{\langle -1\rangle}(x^{-1})}.
\end{equation}
Then, for $k\geq 1$,
\begin{equation*}
\Sigma_k=-\f{1}{k}[x^{-1}]_{\infty}\prod_{j=0}^{k-1}G(x-j).
\end{equation*}
\end{theorem}

Theorem~\ref{bist} implicitly determines $\Sigma_k$ as
a polynomial in the $R_i$'s.
For explicit formulas, it
is convenient to consider separately the graded pieces
of $\Sigma_k$, defined as follows:
let the weight of the monomial $R_{j_1}\ld R_{j_i}$  be $j_1+\ld +j_i$.
For $n\geq 0$, we define
\begin{equation}\label{grade}
\Sigma_{k,2n}=[u^{k+1-2n}]\Sigma_k(R_2u^2,\ld ,R_{k+1}u^{k+1}),
\end{equation}
the sum of all terms of weight $k+1-2n$ in $\Sigma_k$. (From elementary parity
considerations, all other coefficients in $\Sigma_k$ are $0$.)
It is immediate that $\Sigma_{k,0}=R_{k+1}$. An explicit
formula is known for $\Sigma_{k,2}$, and for the statement of this
formula, we introduce polynomials $C_m$ in the $R_i$'s,
where $C_0=1$, $C_1=0$, and
\begin{equation}\label{cexplicit}
C_m=\sum_{\stackrel{j_2,j_3,\ld \geq 0}{2j_2+3j_3+\ld =m}}
\!\!\!\!\!\!\!\! (j_2+j_3+\ld )!\prod_{i\geq 2}\f{((i-1)R_i)^{j_i}}{j_i!},
\;\;\;\; m\geq 2.
\end{equation}
The following explicit formula for $\Sigma_{k,2}$ was
conjectured by Biane~\cite[Conjecture 6.4]{bi},
and proved by \'Sniady~\cite[Theorem 22]{sn}.
\'Sniady's proof was obtained by
finding and then solving an equivalent combinatorial problem.
\begin{theorem}\label{snpf}
For $k\geq 1$,
\begin{equation*}
\Sigma_{k,2}=\tf{1}{24}(k-1)k(k+1)C_{k-1}.
\end{equation*}
\end{theorem}

Note that the R-positivity of $\Sigma_{k,2}$ follows
immediately from Theorem~\ref{snpf}, using~(\ref{cexplicit}).

For $n\geq 2$, only one explicit result is known, given
in the following result for the linear coefficient, due to Biane~\cite{bi}
and Stanley~\cite{st2}.
\begin{theorem}\label{bilin}
For $n\geq 1$, $k\geq 2n-1$, the coefficient of $R_{k+1-2n}$ in $\Sigma_{k,2n}$ is
equal to the number of $k$-cycles $c$ in $\mfS_k$ such that $(1\ld k)c$ has $k-2n$ cycles.
\end{theorem}

Finally, for higher order terms when $n\geq 2$, the following conjecture
of Stanley (private communication) has been communicated to us by Biane.
\begin{conjecture}\label{sthigh}
For $i\geq 1$,
\begin{equation*}
[R_2^i]\Sigma_{2i+3,4}=\tf{1}{540}i(i+1)^3(i+2)^3(i+3)(2i+3).
\end{equation*}
\end{conjecture}

\subsection{Outline of paper}
In this paper, we obtain an explicit
formula for $\Sigma_{k,2n}$, where $k$ and $n$ are arbitrary.
This is our main result, stated in Section 2
as Theorem~\ref{main}. Variants are given also,
as Theorems~\ref{mainmod} and~\ref{maingen}.
These results are a natural generalization of Theorem~\ref{snpf}, since
they give $\Sigma_{k,2n}$ as a polynomial in the $C_m$'s, with
coefficients that are rational polynomials in $k$.
We call such an expression a {\em C-expansion} for $\Sigma_{k,2n}$.
Based on significant amounts of data, we conjecture that
 $\Sigma_{k,2n}$ is {\em C-positive} (all nonzero coefficients
are positive) for all $n\geq 1$, as Conjecture~\ref{Cposall}. This
C-positivity conjecture is stronger than the R-positivity conjecture,
immediately from~(\ref{cexplicit}).

In Section 3, we consider the special cases of our main result for $n=1$
 and $n=2$. For $n=1$, this gives another proof of
Theorem~\ref{snpf}. For $n=2$,
the expression for $\Sigma_{k,4}$ that we obtain, in Theorem~\ref{cor:neq2},
is new. We are able to specialize
this expression to prove Conjecture~\ref{sthigh}.
Also, we are able to prove the C-positivity conjecture for $\Sigma_{k,4}$,
as Corollary~\ref{Cpos}.
Finally, we consider the linear terms in the $R_i$'s, for arbitrary $n$,
and obtain another proof of Theorem~\ref{bilin}.

In general, for $n\geq 3$,
we are not able to prove the R-positivity conjecture nor
the C-positivity conjecture, perhaps because our methods
are not combinatorial.
Instead we apply
Lagrange inversion to ``unwind'' the
compositional inverse in Theorem~\ref{bist}. This
is carried out in Section 4, where we
give the proof of the main result and variants.

\section{The main result}

For the partition $\lam \vdash n$ we denote the {\em monomial} symmetric
function with exponents given by the parts of $\lam$, in
indeterminates $x_1,x_2,\ld$, by $m_{\lam}$. In
this paper, we consider the particular evaluation of the monomial
symmetric function at $x_i=i$, for $i=1,\ld ,k-1$, and $x_i=0$,
for $i\geq k$, and write this as ${\mh}_{\lam}$.
Now let $C(t)=\sum_{m\geq 0}C_mt^m$, so from~(\ref{cexplicit}) we obtain
\begin{equation}\label{Cdef}
C(t)=\f{1}{1-\sum_{i\geq 2}(i-1)R_it^i}.
\end{equation}
Let $D=t\f{d}{dt}$, and $I$ be the identity operator, and define
\begin{equation}
P_m(t)=-\f{1}{m!}C(t)(D+(m-2)I)C(t)\ld (D+I)C(t)DC(t),\;\;\;\; m\geq 1.
\end{equation}
For example, we have
\begin{equation*}\label{P123}
P_1(t)=-C(t),\;\;\;\; P_2(t)=-\f{1}{2}C(t)DC(t),
\end{equation*}
\begin{equation*}
P_3(t)=-\f{1}{6}\left(C(t)^2DC(t)+C(t)(DC(t))^2+C(t)^2D^2C(t)\right).
\end{equation*}
Finally, for a partition $\lam$, we
write $P_{\lam}(t)=\prod_{j=1}^{l(\lam )} P_{\lam_j}(t)$.
We now state our main result.

\begin{theorem}\label{main}
For $n\geq 1$, $k\geq 2n-1$,
\begin{equation*}
\Sigma_{k,2n}=-\f{1}{k}[t^{k+1-2n}]\sum_{\lam\vdash 2n}\mh_{\lam}
\f{P_{\lam}(t)}{C(t)}.
\end{equation*}
\end{theorem}

There is a slight modification
of this result, given below, in which the term corresponding to
the partition with one part is given a simpler (but equivalent)
evaluation.

\begin{theorem}\label{mainmod}
For $n\geq 1$, $k\geq 2n-1$,
\begin{equation*}
\Sigma_{k,2n}=-\f{1}{k}[t^{k+1-2n}]
\left(\f{k-1}{2n}\mh_{2n}P_{2n-1}(t)+\sum_{\stackrel{\lam\vdash 2n}{l(\lam )\geq 2}}\mh_{\lam}
\f{P_{\lam}(t)}{C(t)}\right).
\end{equation*}
\end{theorem}

The following result gives a generating function form of the main result.

\begin{theorem}\label{maingen}
For $n\geq 1$, $k\geq 2n-1$,
\begin{equation*}
\Sigma_{k,2n}=-\f{1}{k}[u^{2n}t^{k+1}]\f{1}{C(t)}\prod_{j=1}^{k-1}
(1+\sum_{i\geq 1}j^iP_i(t)u^it^i),
\end{equation*}
\begin{equation*}
\Sigma_k=-\f{1}{k}[t^{k+1}]\f{1}{C(t)}\prod_{j=1}^{k-1}
(1+\sum_{i\geq 1}j^iP_i(t)t^i).
\end{equation*}
\end{theorem}

Note that, for each $n\geq 1$, these results
give $\Sigma_{k,2n}$ as the coefficient of $t^{k+1-2n}$ in
a polynomial in $C(t)$ and
$$D^iC(t)=\sum_{m\geq 2} m^iC_mt^m,\;\;\;\; i\geq 1.$$
Thus $\Sigma_{k,2n}$ is written as
a polynomial in the $C_m$'s, with coefficients that are
rational in $k$, so our results give C-expansions for
 $\Sigma_{k,2n}$, for $n\geq 1$.

Using the above results, with the help of Maple,
we have determined the C-expansions and the
R-expansions of $\Sigma_{k,2n}$ for
all $k\leq 25$ and $n\geq 1$. The R-expansions are in
complete agreement with those reported in Biane~\cite{bi}
for $k\leq 11$. The C-expansions are given below
for $k\leq 10$:
\begin{eqnarray*}
\Sigma_1-R_2&=&0\\
\Sigma_2-R_3&=&0\\
\Sigma_3-R_4&=&C_2\\
\Sigma_4-R_5&=&\tf{5}{2}C_3\\
\Sigma_5-R_6&=&5\,{C}_{{4}}+8\,{C}_{{2}}\\
\Sigma_6-R_7&=&\tf{35}{4}\,{C}_{{5}}+42\,{C}_{{3}}\\
\Sigma_7-R_8&=& 14\,{C}_{{6}}+{\tf{469}{3}}\,{C}_{{4}}
+{\tf{203}{3}}\,{{C}_{{2}}}^{2} + 180\,{C}_{{2}}\\
\Sigma_8-R_9&=&21\,{C}_{{7}}+{\tf{1869}{4}}\,{C}_{{5}}
+{\tf{819}{2}}\,C_3 C_2+1522\,{C}_{{3}}\\
\Sigma_9-R_{10}&=&30\,{C}_{{8}}+1197\,{C}_{{6}}+
{\tf{963}{2}}\,{{C}_{{3}}}^{2}+1122\,{C}_{{4}}{C}_{{2}}
+81\,{{C}_{{2}}}^{3}+{\tf{26060}{3}}\,{C}_4
+{\tf{17680}{3}}\,{{C}_{{2}}}^{2}\\
&+&8064\,{C}_{{2}}\\
\Sigma_{10}-R_{11}&=&{\tf{165}{4}}\,{C}_{{9}}+{\tf{5467}{2}}\,{C}_{{7}}
+{\tf{4433}{2}}\,{C}_{{4}}{C}_{{3}}
+{\tf{1133}{2}}\,{C}_{{3}}{{C}_{{2}}}^{2}
+{\tf{11033}{4}}
\,{C}_{{5}}{C}_{{2}}+
38225\,{C}_{{5}}\\
&+&52580\,{C}_{{3}}{C}_{{2}}+96624\,{C}_{{3}}
\end{eqnarray*}

Note the form of the data presented above. We have
$$\Sigma_k-\Sigma_{k,0}=\sum_{k\geq 1}\Sigma_{k,2n},$$
where $\Sigma_{k,0}=R_{k+1}$ remains on the lefthandside,
and we can recover the individual $\Sigma_{k,2n}$ on
the righthandside: if the weight of the monomial
 $C_{m_1}\ld C_{m_i}$ is $m_1+\ld +m_i$, then,
from~(\ref{cexplicit}) and~(\ref{grade}),  $\Sigma_{k,2n}$
is the sum of all terms of weight $k+1-2n$.

In the above C-expansions for $k\leq 10$, all nonzero coefficients are positive
rationals, with apparently small denominators. In fact, this is true
for all the data we have computed, up to $k=25$. We do
not have a precise conjecture about the denominators, but conjecture that
the positivity holds for all $k$.

\begin{conjecture}\label{Cposall}
For $n\geq 1$, $k\geq 2n-1$,  $\Sigma_{k,2n}$ is C-positive.
\end{conjecture}

This C-positivity conjecture implies the R-positivity
conjecture, from~(\ref{cexplicit}) (so, our data
also check the R-positivity conjecture for $k\leq 25$). Theorem~\ref{snpf} gives
an immediate proof that Conjecture~\ref{Cposall} holds for $n=1$
 and all $k$. In Corollary~\ref{Cpos}, we are able to prove that
Conjecture~\ref{Cposall} holds for $n=2$ and all $k$. We are not able
to prove the conjecture for any larger value of $n$, though of course
Theorem~\ref{bilin}, together with~(\ref{Cdef}), proves
 that the linear terms are C-positive
for all $n$.

The conjecture does not hold for $n=0$, as described below. We have
 $\Sigma_{k,0}=R_{k+1}$, and it  is straightforward to determine
the C-expansion for the $R_i$'s:
from~(\ref{Cdef}), we obtain
\begin{eqnarray*}
1-\sum_{i\geq 2}(i-1)R_it^i&=&\f{1}{C(t)}\\
&=&\sum_{j_2,j_3,\ld\geq 0}\!\!\!\! (j_2+j_3+\ld )!\prod_{m\geq 2}
\f{(-C_mt^m)^{j_m}}{j_m!},
\end{eqnarray*}
so we conclude that
\begin{equation*}
R_i=\f{1}{i-1}\sum_{\stackrel{j_2,j_3,\ld\geq 0}{2j_2+3j_3+\ld =i}}
\!\!\!\!\!\!\!\! (-1)^{1+j_2+j_3+\ld}(j_2+j_3+\ld )!\prod_{m\geq 2}
\f{C_m^{j_m}}{j_m!},\;\;\;\; i\geq 2.
\end{equation*}
Thus, terms of negative sign appear in the C-expansion of $R_i$,
for $i\geq 4$. This is the reason that we have presented the
data for $k$ up to $10$ with $R_{k+1}$ subtracted on
the lefthandside. This is also the reason that the R-positivity
conjecture does not imply the C-positivity conjecture, so
R-positivity and C-positivity are not equivalent.

\section{Special cases of the main result}

\subsection{Monomial symmetric functions}
To make the expression for $\Sigma_{k,2n}$ that arises from
Theorem~\ref{main} (or Theorem~\ref{mainmod}) explicit, we
need to evaluate the $\mh_{\lam}$, which are monomial symmetric functions
in $1,2,\ld,k-1$. For general results about symmetric functions,
see Macdonald~\cite{ma}.

\begin{proposition}\label{monomial}
For indeterminates $a_i$, $i\geq 1$, let $A(x)=1+\sum_{i\geq 1}a_ix^i$,
and $a_{\lam}=\prod_{j=1}^{l(\lam )}a_{\lam_j}$,
where $\lam=\lam_1\ld\lam_{l(\lam )}$ is a partition. Then
\begin{equation*}
\sum_{\lam\in\mcP}\mh_{\lam}a_{\lam}=
\exp\sum_{j\geq 1}\mh_j\sum_{i\geq 1}\f{(-1)^{i-1}}{i}
[x^j](A(x)-1)^i.
\end{equation*}
\end{proposition}

\noindent
\begin{proof}
We have
\begin{eqnarray*}
\sum_{\lam\in\mcP}m_{\lam}a_{\lam}&=&\prod_{n\geq 1}A(x_n)\\
&=&\exp\sum_{n\geq 1}\log(A(x_n))\\
&=&\exp\sum_{n\geq 1}\sum_{i\geq 1}\f{(-1)^{i-1}}{i}(A(x_n)-1)^i,
\end{eqnarray*}
and the result follows.
\end{proof}

Proposition~\ref{monomial} gives an expression for $\mh_{\lam}$ as
a polynomial in $\mh_i$, $i\geq 1$, by equating coefficients
of $a_{\lam}$. To evaluate the $\mh_i$, $i\geq 1$, we apply
the following result (see, e.g., \cite[I 2, Exercise 11]{ma} for a proof). 

\begin{proposition*}\label{stirling}
For $j\geq 1$,
\begin{equation*}
\mh_j=\sum_{i=1}^jS(j,i)i!{k\choose i+1},
\end{equation*}
where $S(j,i)$, the {\em Stirling numbers of the second kind},
are given by
\begin{equation*}
\sum_{i\geq 0}\sum_{j=0}^iS(j,i)u^i\f{x^j}{j!}
=\exp u(e^x-1).
\end{equation*}
\end{proposition*}

As special cases of this result, we have the following, well-known
sums of integer powers.
\begin{equation}\label{intpowers}
\mh_1=\tf{1}{2}(k-1)k,\;\;\;\; \mh_2=\tf{1}{6}(k-1)k(2k-1),\;\;\;\;
\mh_3=\tf{1}{4}(k-1)^2k^2,
\end{equation}
$$\mh_4=\tf{1}{30}(k-1)k(2k-1)
(3k^2-3k-1).$$

\subsection{The cases $n=1,2$.}
We first consider the
case $n=1$ of Theorem~\ref{mainmod}. This immediately gives
Biane and \'Sniady's C-expansion for $\Sigma_{k,2}$, and hence
another proof of Theorem 1.3, as shown below.
\vspace{.1in}

\noindent
{\bf Proof of Theorem~\ref{snpf}.}
{}From Theorem~\ref{mainmod}, with $n=1$, we obtain
\begin{eqnarray*}
\Sigma_{k,2}&=&-\f{1}{k}[t^{k-1}]
\left( -\tf{1}{2}(k-1)\mh_2C(t)+\mh_{11}C(t)\right)\\
&=&\f{1}{k}\left(\tf{1}{2}(k-1)\mh_2-\mh_{11}\right)[t^{k-1}]C(t).\\
\end{eqnarray*}
But from Proposition~\ref{monomial}, we obtain
$$\mh_{11}=\tf{1}{2}(\mh_1^2-\mh_2),$$
and the result follows from~(\ref{intpowers}), by routine
manipulation.
\hfill$\Box$
\vspace{.1in}

Next we consider the case $n=2$ of Theorem~\ref{mainmod}, to
obtain an explicit C-expansion for $\Sigma_{k,4}$.

\begin{theorem}\label{cor:neq2}
For $k\geq 3$,
\begin{equation*}
\Sigma_{k,4}=\al (k)\!\!\!\!\sum_{\stackrel{i,j,m\geq 0}{i+j+m=k-3}}
\!\!\!\!C_iC_jC_m+
 \be (k)\!\!\!\!\sum_{\stackrel{i,j,m\geq 0}{i+j+m=k-3}}
\!\!\!\! i^2C_iC_jC_m,
\end{equation*}
where
\begin{eqnarray*}
\al (k)&=&-\tf{1}{17280}(k-3)(k-1)^2 k(k+1)(k^2-4k-6),\\
\be (k)&=&\tf{1}{2880}(k-1)k(k+1)(2k^2-3).
\end{eqnarray*}
\end{theorem}

\noindent
\begin{proof}
{}From Theorem~\ref{mainmod}, with $n=2$,
letting $b=\tf{1}{6}(\mh_{31} -\tf{1}{4}(k-1)\mh_4)$, we obtain
\begin{eqnarray*}
\Sigma_{k,4}&=&-\f{1}{k}[t^{k-3}]
(b(C(t)^2DC(t)+C(t)(DC(t))^2+C(t)^2D^2C(t))\\
&+&\tf{1}{4}\mh_{22}C(t)(DC(t))^2-\tf{1}{2}\mh_{211}C(t)^2DC(t)+\mh_{1111}C(t)^3)\\
&=&-\f{1}{k}[t^{k-3}](\mh_{1111}C(t)^3+(b-\tf{1}{2}\mh_{211})C(t)^2DC(t)\\
&+&bC(t)^2D^2C(t)+(b+\tf{1}{4}\mh_{22})C(t)(DC(t))^2)\\
&=&-\f{1}{k}[t^{k-3}](\mh_{1111}C(t)^3+(b-\tf{1}{2}\mh_{211})
\tf{1}{3}DC(t)^3\\
&+&bC(t)^2D^2C(t)+(b+\tf{1}{4}\mh_{22})
(\tf{1}{6}D^2C(t)^3-\tf{1}{2}C(t)^2D^2C(t)))\\
&=&-\f{1}{k}(\mh_{1111}+\tf{1}{3}(k-3)(b-\tf{1}{2}\mh_{211})
+\tf{1}{6}(k-3)^2(b+\tf{1}{4}\mh_{22}))[t^{k-3}]C(t)^3\\
&-&\f{1}{k}(\tf{1}{2}b-\tf{1}{8}\mh_{22})[t^{k-3}]C(t)^2D^2C(t).
\end{eqnarray*}
But from Proposition~\ref{monomial}, we obtain
\begin{eqnarray*}
\mh_{31}&=&\mh_3\mh_1-\mh_4,\\
\mh_{22}&=&\tf{1}{2}(\mh_2^2-\mh_4),\\
\mh_{211}&=&\tf{1}{2}(\mh_2\mh_1^2-2\mh_3\mh_1-\mh_2^2+2\mh_4),\\
\mh_{1111}&=&\tf{1}{24}(\mh_1^4-6\mh_2\mh_1^2+8\mh_3\mh_1+3\mh_2^2-6\mh_4),
\end{eqnarray*}
so from~(\ref{intpowers}), by routine
manipulation, we obtain
\begin{equation}\label{gfS4}
\Sigma_{k,4}=\al(k)[t^{k-3}]C(t)^3+\be(k)[t^{k-3}]C(t)^2D^2C(t),
\end{equation}
where $\al (k)$ and $\be (k)$ are given above. The result follows.
\end{proof}

For monomials in  $R_2,R_3,\ld$ that are pure powers of a single $R_m$,
we have the following form of the above result. 

\begin{corollary}\label{cor:neq3}
For $m\geq 2$, $i\geq 1$,
\begin{eqnarray*}
[R_m^i]\Sigma_{mi+3,4}&=&\tf{1}{34560}(m-1)^i mi(i+1)(i+2)
(mi+2)(mi+3)(mi+4)\\
&\times&(m^3 i^3+2m^2(m+4)i^2+4m(3m+5)i+15m+18).
\end{eqnarray*}
\end{corollary}

\noindent
\begin{proof}
{}From Theorem~\ref{cor:neq2}, we obtain
$$[R_m^i]\Sigma_{mi+3,4}=\al(mi+3)[R_m^it^{mi}]C(t)^3
+\be(mi+3)[R_m^it^{mi}]C(t)^2D^2C(t).$$
Now, setting $R_j=0$ for $j\neq m$, we obtain $C(t)=(1-(m-1)R_mt^m)^{-1}$,
so
$$[R_m^it^{mi}]C(t)^3=(m-1)^i{i+2 \choose 2}.$$
Also, we have
\begin{eqnarray*}
D^2C(t)&=&Dm(m-1)R_mt^m(1-(m-1)R_mt^m)^{-2}\\
&=&Dm\left( (1-(m-1)R_mt^m)^{-2}-(1-(m-1)R_mt^m)^{-1}\right)\\
&=&m^2(m-1)
\left(2R_mt^m(1-(m-1)R_mt^m)^{-3}-R_mt^m(1-(m-1)R_mt^m)^{-2}\right),
\end{eqnarray*}
so
$$[R_m^it^{mi}]C(t)^2D^2C(t)=
(m-1)^im^2\left(2{i+3 \choose 4}-{i+2 \choose 3}\right).$$
The result follows by routine manipulation.
\end{proof}

We now consider the case $m=2$ of Corollary~\ref{cor:neq3}, to obtain
an immediate proof of Stanley's Conjecture~\ref{sthigh}.
\vspace{.1in}

\noindent
{\bf Proof of Conjecture~\ref{sthigh}.}
We set $m=2$ in Corollary~\ref{cor:neq3}. Then the factor that
is cubic in $i$ becomes
$$8i^3+48i^2+88i+48=8(i+1)(i+2)(i+3),$$
and the result follows.
\hfill$\Box$
\vspace{.1in}

As the final result of this section, we are able to use
the explicit C-expansion given in Theorem~\ref{cor:neq2},
to prove the C-positivity of $\Sigma_{k,4}$.

\begin{corollary}\label{Cpos}
$\Sigma_{k,4}$ is $C$-positive for all $k\geq 3$.
\end{corollary}

\noindent
\begin{proof}
Consider $0\leq i\leq j\leq m$, with $i+j+m=k-3$, and let
 $\ga =\vert\Aut (i,j,m)\vert$. Thus when $k=12$, for example,
 $\ga=1$ for $(i,j,m)=(2,3,4)$ or $(0,2,7)$, $\ga =2$ for
 $(i,j,m)=(2,2,5)$ or $(1,4,4)$, and $\ga =6$ for
 $(i,j,m)=(3,3,3)$. Then, from Theorem~\ref{cor:neq2}, we obtain
\begin{equation*}
[C_iC_jC_m]\Sigma_{k,4}=\f{6}{\ga}\al (k)+\f{2}{\ga}(i^2+j^2+m^2)\be (k).
\end{equation*}
Now, the minimum value of $x^2+y^2+z^2$ over the reals, subject to $x+y+z=c$,
for any fixed real $c$, is achieved at $x=y=z=c/3$, so in the above expression
we have $i^2+j^2+m^2\geq \tf{1}{3}(k-3)^2$. But $\be(k)>0$ for $k\geq 3$, so we
obtain
\begin{eqnarray*}
[C_iC_jC_m]\Sigma_{k,4}&\geq &\f{2}{\ga}
\left( 3\al (k)+\tf{1}{3}(k-3)^2\be (k)\right)\\
&=&\tf{1}{8640\ga}(k-3)(k-1)k(k+1)
\left(-3(k-1)(k^2-4k-6)+2(k-3)(2k^2-3)\right)\\
&=& \tf{1}{8640\ga}(k-3)(k-1)k^3(k+1)(k+3)\\
&\geq& 0,
\end{eqnarray*}
for $k\geq 3$, giving the result.
\end{proof}

\subsection{The linear terms.}

We now apply Theorem~\ref{maingen} to evaluate the linear
terms in $\Sigma_k$, and thus obtain another proof
of Theorem~\ref{bilin}.
\vspace{.1in}

\noindent
{\bf Proof of Theorem~\ref{bilin}.}
For $i\geq 1$, let $A^{(i)}(t)$ consist of the terms in $P_i(t)$ that are
linear in the $C_m$'s. Also, let $L_{n,k}=[R_{k+1-2 n}]\Sigma_{k,2n}$.
We apply Theorem~\ref{maingen} to determine $L_{n,k}$.
From~(\ref{Cdef}), we have
\begin{eqnarray*}
L_{n,k}&=&\left[\f{C_{k+1-2n}}{k-2n}\right]\Sigma_{k,2n}=
\left[\f{C_{k+1-2n}}{k-2n}\right]\Sigma_k\\
&=&-\f{1}{k}\left[\f{C_{k+1-2n}}{k-2n}t^{k+1}\right]\f{1}{C(t)}
\prod_{j=1}^{k-1}(1-jt+\sum_{i\geq 1}j^i A^{(i)}(t)t^i)\\
&=&-\f{1}{k}\left[\f{C_{k+1-2n}}{k-2n}t^{k+1}\right]
\f{1}{C(t)}
\left(\prod_{j=1}^{k-1}(1+\sum_{i\geq 1}\f{j^i A^{(i)}(t)t^i}{1-jt})\right)
\prod_{a=1}^{k-1}(1-at)\\
&=&-\f{1}{k}\left[\f{C_{k+1-2n}}{k-2n}t^{k+1}\right]
\left(1-C(t)+\sum_{j=1}^{k-1}\sum_{i\geq 1}\f{j^i A^{(i)}(t)t^i}{1-jt}\right)
\prod_{a=1}^{k-1}(1-at).
\end{eqnarray*}
But
$$A^{(i)}(t)=-\f{1}{i!}(D+(i-2)I)\ld (D+I)DC(t)
=-\sum_{m\geq 2}{-(m-1)\choose i}(-1)^i \f{C_m}{m-1}t^m,\;\;\;\; i\geq 1.$$
Now
let $\f{C_m}{m-1}=x^{m-1}$, $m\geq 2$, which gives
\begin{eqnarray*}
\sum_{i\geq 1}j^i A^{(i)}(t) t^i&=&-\sum_{m\geq 2}
\left((1-jt)^{-(m-1)}-1\right) x^{m-1} t^m\\
&=&-\f{t}{1-\f{xt}{1-jt}}+\f{t}{1-xt},
\end{eqnarray*}
and
$$1-C(t)=-\sum_{m\geq 2}(m-1)x^{m-1}t^m=-\f{t}{(1-xt)^2}
+\f{t}{1-xt}.$$
Thus we obtain
$$L_{n,k}=\f{1}{k}[x^{k-2n}t^{k+1}]
(\f{t}{(1-xt)^2}-\f{t}{1-xt}+\sum_{j=1}^{k-1}
(\f{t}{1-(j+x)t}-\f{t}{(1-jt)(1-xt)}))
\prod_{a=1}^{k-1}(1-at).$$
We now finish the proof using the method of Biane~\cite[Theorem 6.1]{bi}:
Replace $t$ by $t^{-1}$, and multiply by $t^{k}$, to obtain
$$L_{n,k}=\f{1}{k}[x^{k-2n}][t^{-1}]_{\infty}(t)_k
(\f{t}{(t-x)^2}-\f{1}{t-x}+\sum_{j=1}^{k-1}
(\f{1}{t-j-x}-\f{t}{(t-j)(t-x)})),$$
where $(t)_k=t(t-1)\ld (t-k+1)$ is the falling factorial.
Now use the fact that the residue is unchanged if we substitute
 $t+c$ for $t$, where $c$ is independent of $t$. Thus,
substituting $t+j+x$ for $t$ in the first term of the summation
over $j$, and substituting $t+x$ for $t$ in all other terms, we
obtain
\begin{eqnarray*}
L_{n,k}&=&\f{1}{k}[x^{k-2n}]([t](t+x)(t+x)_k-(x)_k+\sum_{j=1}^{k-1}
((x+j)_k-\f{x(x)_k}{x-j}))\\
&=& \f{1}{k}[x^{k-2n}]\sum_{j=0}^{k-1}(x+j)_k
= \f{1}{k}[x^{k-2n}]\sum_{j=0}^{k-1}(x-j)_k,
\end{eqnarray*}
where, for the last equality, we have replaced $x$ by $-x$, and
multiplied by $(-1)^k$. The result now follows,
as shown in Biane~\cite{bi}.
\hfill$\Box$

\section{Lagrange inversion and the proof of the main result}

As a first step, we translate Theorem~\ref{bist} into formal power series,
using the notation
\begin{equation}\label{phiPhi}
\phi(x)=xG(x^{-1}),\;\;\;\;\;\;
\Phi(x,u)=\sum_{i\geq 0}\Phi_i(x)u^i=(1-ux)\phi(x(1-ux)^{-1}),
\end{equation}
where $G(x)$ is defined in~(\ref{FGR}).

\begin{proposition}
  The following two equations hold.
  \vspace{.1in}
  
  \noindent
  1) For $k\geq 1$,
\begin{equation}\label{Six}
\Sigma_k=-\f{1}{k}[x^{k+1}]\prod_{j=0}^{k-1}\Phi(x,j).
\end{equation}

\noindent
2) For $k,n\geq 1$,
\begin{equation}\label{Sixn}
\Sigma_{k,2n}=-\f{1}{k}[u^{2n}x^{k+1}]\prod_{j=0}^{k-1}\Phi(x,ju).
\end{equation}
\end{proposition}

\noindent
\begin{proof}
For~(\ref{Six}), we first replace $x$ by $x^{-1}$ in Theorem~\ref{bist},
to obtain
\begin{equation*}
\Sigma_k=-\f{1}{k}[x^{k+1}]\prod_{j=0}^{k-1}xG(x^{-1}(1-jx)),
\end{equation*}
and the result follows immediately.
\vspace{.1in}

For~(\ref{Sixn}), we let $\vt$ be the substitution
operator $R_i\mpt u^iR_i$, $i\geq 2$. Then, from~(\ref{grade}),
we have 
\begin{equation}\label{subs}
\Sigma_{k,2n}=[u^{k+1-2n}]\vt \Sigma_k .
\end{equation}
Now, from~(\ref{FGR}), we have
\begin{equation*}
\vt F^{\langle -1\rangle}(x)=\f{x}{\vt R(x)}=\f{x}{R(ux)}=\f{1}{u}F^{\langle -1\rangle}(ux),
\end{equation*}
and thus, combining this with~(\ref{FGR}) and~(\ref{phiPhi}), we obtain 
\begin{equation*}
\vt\phi(x)=x\vt G(x^{-1})=\f{x}{\vt F^{\langle -1\rangle}(x)}
=\f{ux}{F^{\langle -1\rangle}(ux)}=\phi(ux),
\end{equation*}
and then
\begin{equation*}
\vt\Phi(x,j)=(1-jx)\phi(ux(1-jx)^{-1})=\Phi(ux,ju^{-1}).
\end{equation*}
Combining this with~(\ref{subs}) and~(\ref{Six}) gives
\begin{equation*}
\Sigma_{k,2n}=-\f{1}{k}[u^{k+1-2n}x^{k+1}]\prod_{j=0}^{k-1}\Phi(ux,ju^{-1})
\end{equation*}
and~(\ref{Sixn}) now follows, by substituting first $x=xu^{-1}$,
and then $u=u^{-1}$.
\end{proof}

Next, we give an expression for the coefficients $\Phi_i$, $i\geq 0$,
defined in~(\ref{phiPhi}).

\begin{proposition}
For $i\geq 0$,
\begin{equation}\label{Phii}
\Phi_i(x)=\f{x}{i!}\left(x^2\f{d}{dx}\right)^i\f{\phi(x)}{x}.
\end{equation}
Note that for $i=0$, this specializes to $\Phi_0(x)=\phi(x)$.
\end{proposition}

\noindent
\begin{proof}
{}From~(\ref{FGR}) and~(\ref{phiPhi}), we have
\begin{equation*}\label{phifps}
\phi(x)=1+\sum_{j\geq 2}\phi_jx^j,
\end{equation*}
where $\phi_j$, $j\geq 2$ are polynomials in the $R_i$'s.
For $i=0$, we have $\Phi_0(x)=\Phi(x,0)=\phi(x)$. For $i\geq 1$,
we have
\begin{eqnarray*}
\Phi_i(x)&=&[u^i]\Phi(x,u)=
[u^i]\left( 1-ux+\sum_{j\geq 2}\phi_jx^j(1-ux)^{1-j}\right)\\
&=&-{1\choose i}x+\sum_{j\geq 2}\phi_j{j+i-2\choose i}x^{j+i}\\
&=&\f{x}{i!}\left(x^2\f{d}{dx}\right) ^i
\left(\f{1}{x}+ \sum_{j\geq 2}\phi_jx^{j-1}\right),
\end{eqnarray*}
and the result follows.
\end{proof}

We make use of the following two, closely related,
versions of Lagrange's Theorem (see, e.g.,~\cite[Section 1.2]{gj}, 
for a proof).

\begin{theorem}\label{lagthm}
Suppose $\psi$ is a formal power series with invertible
constant term. Then the functional equation $s=z\psi(s)$ has
a unique formal power series solution $s=s(z)$. Moreover,

\noindent
1) For a formal Laurent series f and $n\neq 0$, we have
\begin{equation*}
[z^n]f(s)=\f{1}{n}[y^{n-1}]\left(\f{d}{dy}f(y)\right)\psi (y)^n,
\end{equation*}

\noindent
2) For a formal power series $f$, and $n\geq 0$, we have
\begin{equation*}
[z^n]f(s)\f{z}{s}\f{ds}{dz}=[y^n]f(y)\psi (y)^n.
\end{equation*}
\end{theorem}

We consider the functional equation
\begin{equation}\label{lagtw}
w=t\phi(w),
\end{equation}
where $\phi$ is given by~(\ref{phiPhi}). Then from~(\ref{FGR})
and~(\ref{phiPhi}), we have
\begin{equation*}
w=twG(w^{-1})=\f{tw}{F^{\langle -1\rangle}(w)},
\end{equation*}
so $F^{\langle -1\rangle}(w)=t$, and from~(\ref{FGR}) we deduce
that
\begin{equation}\label{lagwt}
t=wR(t).
\end{equation}

We now relate the series $C(t)$ and differential operator $D$ of
Section 2 to the variable $w$.

\begin{proposition}

\begin{equation}\label{wCt}
\f{Dw}{w}=\f{1}{R(t)C(t)}
\end{equation}

\begin{equation}\label{opwt}
w^2\f{d}{dw}=tC(t)D
\end{equation}

\end{proposition}

\begin{proof}
{}From~(\ref{Cdef}) and~(\ref{FGR}), we obtain
\begin{equation*}
C(t)=\f{1}{-tD\f{R(t)}{t}}.
\end{equation*}
But
\begin{equation*}
\f{Dw}{w}=-wD\f{1}{w}=-\f{t}{R(t)}D\f{R(t)}{t},
\end{equation*}
from~(\ref{lagwt}), and result~(\ref{wCt}) follows.
\vspace{.1in}

Now,~(\ref{wCt}) gives the operator identity
\begin{equation*}
w\f{d}{dw}=R(t)C(t)D,
\end{equation*}
and multiplying by $w$ and using~(\ref{lagwt}), we obtain result~(\ref{opwt}).
\end{proof}

\noindent
{\bf Proof of Theorem~\ref{main}.}
For a partition $\lam$, let $\Phi _{\lam}(x)=\prod_{j=1}^{l(\lam)}\Phi _{\lam_j}(x)$. Then
from~(\ref{Sixn}) and~(\ref{Phii}), we have
\begin{eqnarray*}
\Sigma_{k,2n}&=&-\f{1}{k}[x^{k+1}]\sum_{\lam\vdash 2n}\mh_{\lam} \Phi _{\lam}(x)
\phi (x)^{k-l(\lam)}\\
&=&-\f{1}{k}[x^{k+1}]\sum_{\lam\vdash 2n}\mh_{\lam}
 \f{\Phi _{\lam}(x)}{\phi (x)^{l(\lam )+1}}\phi (x)^{k+1}\\
&=&-\f{1}{k}[t^{k+1}]\sum_{\lam\vdash 2n}\mh_{\lam}
 \f{1}{R(t)C(t)}\f{\Phi _{\lam}(w)}{\phi (w)^{l(\lam )+1}},
\end{eqnarray*}
where the last equality follows from Theorem~\ref{lagthm}.2 and~(\ref{wCt}).
But, from~(\ref{Phii}),~(\ref{lagtw}) and~(\ref{opwt}), for $i\geq 1$ we have
\begin{eqnarray*}
\f{\Phi _i(w)}{\phi (w)}&=&\f{1}{i!}\f{w}{\phi (w)}(w^2\f{d}{dw})^i\f{\phi (w)}{w}\\
&=&\f{t}{i!}(tC(t)D)^{i-1}tC(t)D\f{1}{t}\\
&=&-\f{t}{i!}(tC(t)D)^{i-1}C(t).
\end{eqnarray*}

Finally, we prove by induction on $i\geq 1$ that
\begin{equation*}
-\f{1}{i!}(tC(t)D)^{i-1}C(t)=t^{i-1}P_i(t),
\end{equation*}
where $P_i(t)$ is defined in Section 2. The result is clearly
true for $i=1$. For the induction step, we have
\begin{eqnarray*}
-\f{1}{(i+1)!}(tC(t)D)^iC(t)&=&\f{1}{i+1}tC(t)Dt^{i-1}P_i(t)\\
&=&\f{1}{i+1}\left( t^iC(t)D+(i-1)t^iC(t)I\right) P_i(t)\\
&=&t^iP_{i+1}(t),
\end{eqnarray*}
as required. Together, these results give
\begin{equation*}
\f{\Phi _i(w)}{\phi (w)}=t^iP_i(t),
\end{equation*}
so
\begin{equation*}
\f{\Phi _{\lam}(w)}{\phi (w)^{l(\lam )+1}}=t^{2n}\f{P_{\lam}(t)}{\phi (w)},
\end{equation*}
since $\lam\vdash 2n$, and the result follows from \eqref{lagtw} and \eqref{lagwt}.\hfill$\Box$
\vspace{0.1in}

\noindent
{\bf Proof of Theorem~\ref{mainmod}.}
In the proof of Theorem~\ref{main},
the term in $\Sigma_{k,2n}$ corresponding to the partition with the single part $2n$ can
be treated in the following modified way. We obtain
\begin{eqnarray*}
-\f{1}{k}[x^{k+1}]\mh_{2n} \Phi _{2n}(x)\phi (x)^{k-1}&=&
-\f{1}{k}[x^{k-2}]\mh_{2n} x^{-3}\Phi _{2n}(x)\phi (x)^{k-1}\\
&=&-\f{1}{k}[x^{k-2}]\mh_{2n}x^{-3}\f{x}{(2n)!}x^2\f{d}{dx}
\left(x^2\f{d}{dx}\right)^{2n-1}\f{\phi (x)}{x}\\
&=&-\f{k-1}{k}[t^{k-1}]\mh_{2n}\f{1}{(2n)!}\left(w^2\f{d}{dw}\right)^{2n-1}\f{\phi (w)}{w},
\end{eqnarray*}
from Theorem~\ref{lagthm}.1, and the  result follows as in the above proof of
Theorem~\ref{main}.\hfill$\Box$

\section*{Acknowledgements}

This work was supported by a Discovery Grant from NSERC (IG), a
Postgraduate Scholarship from NSERC (AR), and an Ontario Graduate Scholarship
in Science and Technology (AR). We would like to thank
P. Biane, A. Okounkov, P. \'Sniady and R. Stanley for helpful
comments on an earlier draft.

\bibliographystyle{hplain}
\bibliography{kerovpolys}

\begin{thebibliography}{10}

\bibitem{bi1}
P.~Biane.
\newblock Representations of the symmetric groups and free probability.
\newblock {\em Advances in Mathematics}, 138:126--181, 1998.

\bibitem{bi2}
P.~Biane.
\newblock Free cumulants and representations of large symmetric groups.
\newblock {\em Proceedings of the XIIIth International Congress of Mathematical
  Physics, London,Int. Press}, pages 321--326, 2000.

\bibitem{bi}
P.~Biane.
\newblock Characters of symmetric groups and free cumulants.
\newblock {\em Asymptotic Combinatorics with Applications to Mathematical
  Physics, A. Vershik (Ed.), Springer Lecture Notes in Mathematics,},
  1815:185--200, 2003.

\bibitem{cgs}
S.~Corteel, A.~Goupil, and G.~Schaeffer.
\newblock Content evaluation and class symmetric functions.
\newblock preprint, 2004.

\bibitem{fjr}
A.~Frumkin, G.~James, and Y.~Roichman.
\newblock On trees and characters.
\newblock preprint, 2001.

\bibitem{gj}
I.P. Goulden and D.M. Jackson.
\newblock {\em Combinatorial Enumeration}.
\newblock Wiley-Interscience, New York, 1983 (Dover Reprint, 2004).

\bibitem{io}
V.~Ivanov and G.~Olshanski.
\newblock Kerov's central limit theorem for the plancherel measure on young
  diagrams.
\newblock {\em Symmetric functions 2001: Surveys of developments and
  perspectives, S. Fomin (Ed.), NATO Science series II. Mathematics, Physics
  and Chemistry}, 74:93--151, 2002, arXiv:math.CO/0304010.

\bibitem{ka}
J.~Katriel.
\newblock Explicit expressions for the central characters of the symmetric
  group.
\newblock {\em Discrete Applied Math.}, 67:149--156, 1996.

\bibitem{ke}
S.~Kerov.
\newblock Gaussian limit for the plancherel measure of the symmetric group.
\newblock {\em Comptes Rendus Acad. Sci. Paris, S{\' e}rie}, 316:303--308,
  1993.

\bibitem{ma}
I.G. Macdonald.
\newblock {\em Symmetric Functions and Hall Polynomials}.
\newblock Oxford University Press, Oxford, 2 edition, 1995.

\bibitem{sn}
P.~\'Sniady.
\newblock Asymptotics of characters of symmetric groups and free probability.
\newblock preprint, 2003, arXiv:math.CO/0304275.

\bibitem{st2}
R.P. Stanley.
\newblock Kerov's character polynomial and irreducible symmetric group
  characters of rectangular shape.
\newblock Transparencies from a talk at CMS meeting, Quebec City, 2002.

\end{thebibliography}

\end{document}